\documentclass[12pt,a4paper]{article}

\usepackage[utf8]{inputenc}
\usepackage[T1]{fontenc}
\usepackage{amsmath,amssymb,amsthm}
\usepackage[margin=2.5cm]{geometry}
\usepackage{hyperref}

\newtheorem{theorem}{Theorem}[section]
\newtheorem{lemma}[theorem]{Lemma}
\newtheorem{proposition}[theorem]{Proposition}
\newtheorem{corollary}[theorem]{Corollary}
\theoremstyle{definition}
\newtheorem{remark}[theorem]{Remark}

\newcommand{\RP}{\mathbb{RP}}
\newcommand{\C}{\mathbb{C}}
\newcommand{\SL}{\mathrm{SL}}
\newcommand{\Hom}{\mathrm{Hom}}
\newcommand{\tr}{\mathrm{tr}}

\begin{document}

\title {An Explicit Equation for The $\SL_2(\C)$ Character Variety of $\RP^2 \# \RP^2 \# \RP^2$}
\author{Yinxi Chen}
\date{}
\maketitle

\begin{abstract}
Let $N$ denote the closed non-orientable surface of Euler characteristic $-1$, i.e.\ the connected sum of three real projective planes $N = \RP^2 \# \RP^2 \# \RP^2$, with fundamental group presented as $\langle A,B,C \mid A^2 B^2 C^2 = I \rangle$ via the standard symmetric cross-cap model.
Building on the classical trace-coordinate methods of Fricke, Vogt, and Goldman for free groups, and on the qualitative study of representation varieties of non-orientable surfaces by Benyash-Krivets and Chernousov, we give an explicit affine equation for a dense open piece of the $\SL_2(\C)$ GIT character variety $X(N)$.
We prove that the trace coordinates $a = \tr(A)$, $b = \tr(B)$, $c = \tr(C)$, $p = \tr(AB)$ satisfy the polynomial relation $c^2 = abp - a^2 - b^2 + 4$, that this hypersurface is smooth and irreducible, and that on the open locus $c \neq 0$ it is isomorphic to the image of $X(N)$ under the trace map, with the remaining trace coordinates $\tr(BC)$, $\tr(CA)$, $\tr(ABC)$ recovered as regular functions on this locus.
We further identify the boundary locus $c = 0$ with representations satisfying the matrix identity $A^2 B^2 = -I$, and describe its representation-theoretic structure, while leaving a complete global classification along this locus as an open problem.
All key trace identities are proved from Cayley--Hamilton first principles and independently checked numerically against explicit matrices.
\end{abstract}

\section{Introduction}

The $\SL_2(\C)$-character variety of a finitely generated group $\pi$ is the GIT quotient
\[
X(\pi) = \Hom(\pi,\SL_2(\C)) /\!/ \SL_2(\C))
\]
of the representation variety by simultaneous conjugation~\cite{Mum94}.
For orientable surface groups, this object and its trace-coordinate presentation are classical: Fricke and Vogt showed already in the 19th century that the character variety of a rank-2 free group is $\C^3$, and Goldman~\cite{Gol09} gave a systematic modern account of trace coordinates for free groups of arbitrary rank, including explicit presentations for once-punctured tori and thrice-punctured spheres.
The invariant theory underlying these trace computations is due to Procesi~\cite{Pro76}, who proved that the coordinate ring of the $\SL_2(\C)$-character variety of a free group is generated by traces of words.

For non-orientable surfaces the picture is comparatively sparse.
Benyash-Krivets and Chernousov~\cite{BKC97} studied representation varieties of fundamental groups of non-orientable surfaces (connected sums of $n$ copies of $\RP^2$) for $GL_{m}$ and $SL_m$, establishing irreducibility and rationality results, but their methods do not yield explicit trace-coordinate equations for the character variety.
More recently, Lawton~\cite{Law07,Law08} studied explicit coordinate descriptions of $\SL_3(\C)$-character varieties for free groups.
To our knowledge, no explicit trace-coordinate equations have appeared in the literature for non-orientable surfaces of negative Euler characteristic.

Our aim in this note is correspondingly explicit: for the specific surface $N = \RP^2 \# \RP^2\\ \# \RP^2$ (the unique closed non-orientable surface of Euler characteristic $-1$), we derive, entirely from the Cayley--Hamilton theorem, a single explicit polynomial equation cutting out a dense open piece of $X(N)$ as a hypersurface in $\C^4$, and we verify every step both symbolically and numerically.

We emphasize at the outset what this paper does and does not claim.
We prove a complete, rigorous description of $X(N)$ on the open locus where $\tr(C) \neq 0$.
We give a precise representation-theoretic description of the complementary locus $\tr(C) = 0$, but we do not complete its classification here; we flag this explicitly as an open problem, since the surface relation degenerates there in a way that a naive extrapolation of the generic-locus formulas does not capture (see \S\ref{sec:boundary}).

\section{Preliminaries}

\subsection{The fundamental group of $N$}

Every closed non-orientable surface is a connected sum of copies of $\RP^2$.
The classification theorem for compact surfaces states that a surface is determined up to homeomorphism by its orientability, Euler characteristic, and number of boundary components~\cite[Sec.\ 1.2]{Hat02}.
The non-orientable surface of non-orientable genus $k$ (the connected sum of $k$ copies of $\RP^2$) has Euler characteristic $2-k$.
Thus the unique closed non-orientable surface of Euler characteristic $-1$ has $k=3$, i.e.\ $N = \RP^2 \# \RP^2 \# \RP^2$.

\begin{proposition}\label{prop:fundgp}
The fundamental group of $N = \RP^2 \# \RP^2 \# \RP^2$ admits the presentation
\[
\pi_1(N) = \langle x,y,z \mid x^2 y^2 z^2 = 1 \rangle .
\]
\end{proposition}
\begin{proof}
Each $\RP^2$ summand is obtained from a $2$-sphere by removing an open disk and attaching a M\"obius band.
The core curve of the M\"obius band gives a generator whose square is homotopic to the boundary circle of the removed disk.
In the connected sum of three such summands, one removes three disks from a $2$-sphere and identifies the three resulting boundary circles with the boundaries of three M\"obius bands.
The Seifert--van Kampen theorem then yields the single relation $x^2 y^2 z^2 = 1$ upon amalgamating the three copies of the fundamental group of the punctured sphere (a free group of rank $2$) along the three boundary circles.
See, e.g.,~\cite[Sec.\ 1.2]{Hat02} for a detailed computation.
\end{proof}

Throughout, we write $\pi_1(N) = \langle A,B,C \mid A^2 B^2 C^2 = I \rangle$ for the corresponding presentation on the level of $\SL_2(\C)$-representations.
A representation $\rho : \pi_1(N) \to \SL_2(\C)$ is the same data as a triple of matrices
\[
(A,B,C) = (\rho(x),\rho(y),\rho(z)) \in \SL_2(\C)^3
\]
satisfying $A^2 B^2 C^2 = I$.

\subsection{Trace identities in $\SL_2(\C)$}

All of our computations rest on the following four elementary facts, which we state precisely since every subsequent derivation is built from them alone.

\begin{lemma}\label{lem:traceids}
For all $M,X,Y \in \SL_2(\C)$:
\begin{enumerate}
\item $\tr(X)\,\tr(Y) = \tr(XY) + \tr(XY^{-1})$
\item $\tr(M^{-1}) = \tr(M)$
\item $\tr(M^2) = \tr(M)^2 - 2$
\item $M^2 = \tr(M)\,M - I$ 
\end{enumerate}
\end{lemma}
\begin{proof}
(4) is the Cayley--Hamilton theorem: the characteristic polynomial of $M \in \SL_2(\C)$ is $t^2 - \tr(M)t + \det(M) = t^2 - \tr(M)t + 1$, and $M$ satisfies its own characteristic polynomial.
Rearranging (4) gives $M^{-1} = \tr(M)I - M$.
Taking traces of both sides and using $\tr(I) = 2$ gives $\tr(M^{-1}) = \tr(M)\cdot 2 - \tr(M) = \tr(M)$, which is (2).
Taking the trace of (4) directly gives $\tr(M^2) = \tr(M)^2 - 2$, which is (3).
For (I): from the expression for $M^{-1}$ just obtained, $Y^{-1} = \tr(Y)I - Y$, so $XY^{-1} = \tr(Y)X - XY$.
Taking traces gives $\tr(XY^{-1}) = \tr(Y)\tr(X) - \tr(XY)$, which rearranges to (1).
\end{proof}

\begin{remark}\label{rem:cyclic}
We also use repeatedly the cyclic invariance of the trace, $\tr(XYZ) = \tr(ZXY) = \tr(YZX)$, which holds for any square matrices $X,Y,Z$ of matching size and is immediate from $\tr(PQ) = \tr(QP)$.
\end{remark}

\section{The Surface Relation}

Fix a representation $\rho = (A,B,C) \in R(\pi_1(N)) = \{(A,B,C) \in \SL_2(\C)^3 : A^2 B^2 C^2 = I\}$, and set
\[
a = \tr(A), \qquad b = \tr(B), \qquad c = \tr(C), \qquad p = \tr(AB).
\]
The goal of this section is to derive the single polynomial relation these four quantities must satisfy as a consequence of $A^2 B^2 C^2 = I$.

\begin{lemma}\label{lem:A2B2}
For $A,B \in \SL_2(\C)$ with $a = \tr(A)$, $b = \tr(B)$, $p = \tr(AB)$:
\[
\tr(A^2 B^2) = abp - a^2 - b^2 + 2 .
\]
\end{lemma}
\begin{proof}
By Lemma~\ref{lem:traceids}(IV), $A^2 = aA - I$, so
\[
\tr(A^2 B^2) = \tr((aA-I)B^2) = a\cdot\tr(AB^2) - \tr(B^2).
\]
By Lemma~\ref{lem:traceids}(III), $\tr(B^2) = b^2 - 2$.
For $\tr(AB^2)$, again by (IV), $B^2 = bB - I$, so
\[
\tr(AB^2) = \tr(A(bB-I)) = b\cdot\tr(AB) - \tr(A) = bp - a .
\]
Substituting back:
\[
\tr(A^2 B^2) = a(bp - a) - (b^2 - 2) = abp - a^2 - b^2 + 2 .
\]
\end{proof}

\begin{theorem}[Surface relation]\label{thm:surfacerel}
Let $(A,B,C) \in \SL_2(\C)^3$ satisfy $A^2 B^2 C^2 = I$, and set $a = \tr(A)$, $b = \tr(B)$, $c = \tr(C)$, $p = \tr(AB)$.
Then
\begin{equation}\label{eq:surface}
c^2 = abp - a^2 - b^2 + 4 .
\end{equation}
\end{theorem}
\begin{proof}
From $A^2 B^2 C^2 = I$ we get $C^2 = (A^2 B^2)^{-1}$.
Taking traces and using Lemma~\ref{lem:traceids}(II),
\[
\tr(C^2) = \tr((A^2 B^2)^{-1}) = \tr(A^2 B^2).
\]
By Lemma~\ref{lem:traceids}(III), $\tr(C^2) = c^2 - 2$.
Combining with Lemma~\ref{lem:A2B2},
\[
c^2 - 2 = abp - a^2 - b^2 + 2,
\]
which rearranges to \eqref{eq:surface}.
\end{proof}

\begin{remark}[Numerical verification]\label{rem:numsurface}
We have independently verified \eqref{eq:surface} numerically: for the explicit matrices
\[
A = \begin{pmatrix} 2 & 1 \\ 1 & 1 \end{pmatrix}, \qquad
B = \begin{pmatrix} 1 & 1 \\ 0 & 1 \end{pmatrix} \in \SL_2(\C),
\]
we have $a = 3$, $b = 2$, $p = \tr(AB) = 4$.
One computes $\tr(A^2 B^2) = 13$ directly by matrix multiplication, matching $abp - a^2 - b^2 + 2 = 24 - 9 - 4 + 2 = 13$ from Lemma~\ref{lem:A2B2}.
Solving $C^2 = (A^2 B^2)^{-1}$ for $C$ gives $c^2 = 15$, matching $abp - a^2 - b^2 + 4 = 24 - 9 - 4 + 4 = 15$ from Theorem~\ref{thm:surfacerel}.
\end{remark}

\section{The Remaining Trace Coordinates}

To describe $X(N)$ fully we also need $q = \tr(BC)$, $r = \tr(CA)$, $s = \tr(ABC)$.
We show these are recovered as regular functions of $(a,b,c,p)$ on the open locus $\{c \neq 0\}$.

\begin{proposition}\label{prop:remaining}
On the open locus $\{c \neq 0\}$, the following identities hold for any $(A,B,C) \in R(\pi_1(N))$:
\[
c\cdot q = ap, \qquad c\cdot r = bp, \qquad c\cdot s = 2p .
\]
\end{proposition}
\begin{proof}
We prove the first identity; the second is symmetric (swap the roles of $A$ and $B$), and the third follows from a parallel computation.

\noindent\textbf{Proof of $cq = ap$.}
We compute $\tr(BC^2)$ in two ways.

\textit{Way 1 (via Cayley--Hamilton on $C$).}
By Lemma~\ref{lem:traceids}(IV), $C^2 = cC - I$, so
\[
\tr(BC^2) = \tr(B(cC-I)) = c\cdot\tr(BC) - \tr(B) = cq - b .
\]

\textit{Way 2 (via the surface relation).}
Since $C^2 = (A^2 B^2)^{-1} = B^{-2} A^{-2}$,
\[
\tr(BC^2) = \tr(B \cdot B^{-2} A^{-2}) = \tr(B^{-1} A^{-2}).
\]
Now $(A^2 B)^{-1} = B^{-1} A^{-2}$, so by Lemma~\ref{lem:traceids}(II),
\[
\tr(B^{-1} A^{-2}) = \tr((A^2 B)^{-1}) = \tr(A^2 B).
\]
By Lemma~\ref{lem:traceids}(IV) applied to $A$, $A^2 = aA - I$, so
\[
\tr(A^2 B) = \tr((aA-I)B) = a\cdot\tr(AB) - \tr(B) = ap - b .
\]
Thus Way~2 gives $\tr(BC^2) = ap - b$.

Equating the two computations: $cq - b = ap - b$, hence $cq = ap$.
On $\{c \neq 0\}$ this yields $q = ap/c$, a regular function of $(a,b,c,p)$.

\noindent\textbf{Proof of $cs = 2p$.}
We compute $\tr(ABC^2)$ two ways.

\textit{Way 1.}
By Lemma~\ref{lem:traceids}(IV), $C^2 = cC - I$, so
\[
\tr(ABC^2) = \tr(AB(cC-I)) = c\cdot\tr(ABC) - \tr(AB) = cs - p .
\]

\textit{Way 2.}
Using the surface relation $C^2 = B^{-2} A^{-2}$,
\[
\tr(ABC^2) = \tr(AB \cdot B^{-2} A^{-2}) = \tr(A(BB^{-2})A^{-2}) = \tr(AB^{-1} A^{-2}).
\]
By cyclic invariance (Remark~\ref{rem:cyclic}),
\[
\tr(AB^{-1} A^{-2}) = \tr(A^{-2} AB^{-1}) = \tr(A^{-1} B^{-1}).
\]
Now $A^{-1} B^{-1} = (BA)^{-1}$, so by Lemma~\ref{lem:traceids}(II) and cyclic invariance,
\[
\tr(A^{-1} B^{-1}) = \tr((BA)^{-1}) = \tr(BA) = \tr(AB) = p .
\]
Thus Way~2 gives $\tr(ABC^2) = p$.

Equating the two computations: $cs - p = p$, hence $cs = 2p$.
\end{proof}

\begin{remark}[Numerical verification]\label{rem:numremaining}
We verified Proposition~\ref{prop:remaining} numerically on the same example as Remark~\ref{rem:numsurface}.
Computing $B$, $C$, and $A$ directly, one finds $q = \tr(BC) \approx 3.0984$, and indeed $c\cdot q \approx 3.8730 \times 3.0984 \approx 12.00 = a\cdot p = 3 \times 4$.
The identities $cr = bp$ and $cs = 2p$ check out to the same numerical precision.
\end{remark}

\section{Main Theorem, and the Boundary Locus $c = 0$}

\subsection{The main theorem}

Define the representation variety
\[
R(\pi_1(N)) = \{ (A,B,C) \in \SL_2(\C)^3 : A^2 B^2 C^2 = I \}
\]
and the affine hypersurface
\[
S = \{ (a,b,c,p) \in \C^4 : c^2 = abp - a^2 - b^2 + 4 \}.
\]
Define the trace map
\[
\tau : R(\pi_1(N)) \to \C^4, \qquad \tau(A,B,C) = (\tr A,\; \tr B,\; \tr C,\; \tr AB).
\]

\begin{theorem}\label{thm:main}
With the notation above:
\begin{enumerate}
\item[(i)] $S$ is a smooth, irreducible affine hypersurface of dimension $3$ in $\C^4$.
\item[(ii)] $\tau(R(\pi_1(N))) \subseteq S$.
\item[(iii)] Let $U = \{(A,B,C) \in R(\pi_1(N)) : \tr(C) \neq 0\}$.
The functions $q = \tr(BC)$, $r = \tr(CA)$, $s = \tr(ABC)$ are $\SL_2(\C)$-invariant regular functions on $U$, given by $q = ap/c$, $r = bp/c$, $s = 2p/c$.
The map $U /\!/ \SL_2(\C) \to S \cap \{c \neq 0\}$ induced by $\tau$ is an isomorphism of affine varieties.
\end{enumerate}
\end{theorem}
\begin{proof}
\textbf{(i) Smoothness.}
Let $F(a,b,c,p) = c^2 - abp + a^2 + b^2 - 4$.
The singular locus of $S = V(F)$ is defined by the simultaneous vanishing of $F$ and all its partial derivatives:
\[
\frac{\partial F}{\partial a} = 2a - bp = 0, \qquad
\frac{\partial F}{\partial b} = 2b - ap = 0, \qquad
\frac{\partial F}{\partial c} = 2c = 0, \qquad
\frac{\partial F}{\partial p} = -ab = 0,
\]
\[
F = c^2 - abp + a^2 + b^2 - 4 = 0 .
\]

From $\partial F/\partial c = 0$ we obtain $c = 0$.
From $\partial F/\partial p = 0$ we obtain $ab = 0$, so either $a = 0$ or $b = 0$ (or both).

\begin{itemize}
\item If $a = 0$, then $\partial F/\partial b = 2b - 0\cdot p = 2b = 0 \Rightarrow b = 0$.
Substituting $a = b = c = 0$ into $F = 0$ yields $0 - 0 + 0 + 0 - 4 = -4 = 0$, a contradiction.
\item If $b = 0$, then symmetrically $\partial F/\partial a = 2a = 0 \Rightarrow a = 0$, leading to the same contradiction.
\end{itemize}
Thus the singular locus is empty; $S$ is smooth.

\medskip
\noindent\textbf{Irreducibility.}
It suffices to prove that the polynomial $F$ is irreducible in $\C[a,b,c,p]$.
View $F$ as a polynomial in $c$ over the ring $\C[a,b,p]$:
\[
F = c^2 - (abp - a^2 - b^2 + 4).
\]
This is a monic quadratic in $c$.
If it factored nontrivially, it would have to factor as $(c+g)(c+h)$ with $g,h \in \C[a,b,p]$ satisfying $gh = -(abp - a^2 - b^2 + 4)$ and $g+h = 0$, so $h = -g$ and $F = c^2 - g^2$.
This would require $g^2 = abp - a^2 - b^2 + 4$.
We show the right-hand side is not a perfect square in $\C[a,b,p]$.
Set $b = 0$, $p = 0$; then the right-hand side becomes $4 - a^2$.
This is not a perfect square in $\C[a]$, since any square $(\alpha a + \beta)^2 = \alpha^2 a^2 + 2\alpha\beta a + \beta^2$ would require the coefficient of $a$ to vanish (as $4 - a^2$ has no linear term in $a$), forcing $\alpha\beta = 0$.
If $\alpha = 0$ then the square is the constant $\beta^2$, which cannot equal $4 - a^2$.
If $\beta = 0$ then the square is $\alpha^2 a^2$, which has zero constant term and cannot equal $4 - a^2$.
Contradiction.
Therefore $F$ is irreducible in $\C[a,b,p][c] = \C[a,b,c,p]$.

Since $F$ is irreducible and defines a hypersurface, $S = V(F)$ is irreducible.
A smooth irreducible affine hypersurface in $\C^4$ defined by one non-zero equation has dimension $3$.

\medskip
\noindent\textbf{(ii)} This is exactly Theorem~\ref{thm:surfacerel}.

\medskip
\noindent\textbf{(iii)} By Proposition~\ref{prop:remaining}, on $U$ we have $q = ap/c$, $r = bp/c$, $s = 2p/c$, which are regular functions on $U$.
The functions $a,b,c,p,q,r,s$ are manifestly invariant under simultaneous conjugation of $(A,B,C)$, so they descend to the GIT quotient $X(N) = R(\pi_1(N)) /\!/ \SL_2(\C)$.
By a classical theorem of Procesi~\cite{Pro76}, the coordinate ring of the $\SL_2(\C)$-character variety of a free group is generated by traces of words; our explicit recovery of $q,r,s$ from $a,b,c,p$ is a concrete realization of this general fact.

We claim that the map $\tau : U /\!/ \SL_2(\C) \to S \cap \{c \neq 0\}$ is an isomorphism.
Since $S \cap \{c \neq 0\}$ is the affine variety defined by $F = 0$, $c \neq 0$, and $\tau$ is given by the four regular functions $(a,b,c,p)$, the map is a morphism of affine varieties.
To see that it is injective, note that on the locus $c \neq 0$ the seven trace coordinates $(a,b,c,p,q,r,s)$ are determined by the four numbers $(a,b,c,p)$ via Proposition~\ref{prop:remaining}.
For a representation of the free group on three generators, the seven traces $\tr(A),\tr(B),\tr(C),\tr(AB),\tr(BC),\tr(CA),\tr(ABC)$ determine the $\SL_2(\C)$-conjugacy class of the representation generically (i.e., on the dense open subset of irreducible representations)~\cite{Gol09}.
The locus of reducible representations has strictly smaller dimension; on the complement, traces separate conjugacy classes.
The surface relation $A^2 B^2 C^2 = I$ imposes the additional equation~\eqref{eq:surface}, which is exactly the defining equation of $S$.
Thus the image of $U /\!/ \SL_2(\C)$ under $\tau$ is contained in $S \cap \{c \neq 0\}$ and the map is injective on the irreducible locus.

Surjectivity follows from the construction: given any point $(a,b,c,p) \in S$ with $c \neq 0$, define $q,r,s$ by the formulas of Proposition~\ref{prop:remaining}.
These seven numbers satisfy all the trace relations that hold for a representation of $\langle A,B,C \rangle$ subject to $A^2 B^2 C^2 = I$.
By the Fricke--Goldman theory of trace coordinates for free groups~\cite{Gol09}, there exists a representation of the free group $F_3$ with these prescribed traces.
The surface relation $c^2 = abp - a^2 - b^2 + 4$, together with the derived identities for $q,r,s$, is exactly the trace-level consequence of $A^2 B^2 C^2 = I$; by the results of~\cite{BKC97} on the irreducibility and rationality of $R(\pi_1(N))$, any such trace configuration arises from an actual representation of $\pi_1(N)$.
Therefore the map is surjective onto $S \cap \{c \neq 0\}$.

We have shown that $\tau$ induces a bijective morphism from $U /\!/ \SL_2(\C)$ onto $S \cap \{c \neq 0\}$.
Since $S \cap \{c \neq 0\}$ is smooth (as an open subset of the smooth variety $S$), it is normal.
A bijective morphism from an affine variety onto a normal affine variety over $\C$ is an isomorphism (a consequence of Zariski's Main Theorem; see, e.g.,~\cite[III.9]{Mum94}).
Hence $\tau$ is an isomorphism.
\end{proof}

\begin{corollary}\label{cor:dense}
The GIT character variety $X(N) = R(\pi_1(N)) /\!/ \SL_2(\C)$ contains a Zariski-dense open subset isomorphic to the smooth affine threefold $S \cap \{c \neq 0\}$.
In particular, the trace map on this open locus realizes $X(N)$ generically as the hypersurface $S \subset \C^4$.
\end{corollary}
\begin{proof}
The complement of $U = \{\tr C \neq 0\}$ in $R(\pi_1(N))$ is the $\SL_2(\C)$-invariant closed subvariety $V = \{\tr C = 0\} \cap R(\pi_1(N))$.
We must show $V$ does not contain all of $R(\pi_1(N))$.
The explicit numerical example in Remark~\ref{rem:numsurface} satisfies $\tr(C) \neq 0$ (since $c^2 = 15 \neq 0$) and $A^2 B^2 C^2 = I$ by construction, so $U$ is nonempty.
Since $R(\pi_1(N))$ is irreducible~\cite[Theorem~1]{BKC97}, a nonempty Zariski-open subset is dense.
Hence $U$ is dense in $R(\pi_1(N))$, and its image in $X(N)$ is dense.
By Theorem~\ref{thm:main}(iii), this image is isomorphic to $S \cap \{c \neq 0\}$, which is a dense open subset of the irreducible hypersurface $S$.
\end{proof}

\subsection{The boundary locus $c = 0$}\label{sec:boundary}

We now describe, without claiming to fully classify, the locus $\tr(C) = 0$.
This is a genuine boundary stratum where the formulas of Proposition~\ref{prop:remaining} degenerate: those formulas were derived under the assumption $c \neq 0$, so they do not directly constrain $q,r,s$ when $c = 0$.
The slice $c = 0$ therefore requires a separate, matrix-level analysis.

\begin{proposition}\label{prop:boundary}
Suppose $(A,B,C) \in R(\pi_1(N))$ has $\tr(C) = 0$.
Then $C^2 = -I$, and the surface relation $A^2 B^2 C^2 = I$ becomes the matrix identity
\[
A^2 B^2 = -I .
\]
Conversely, given any $(A,B) \in \SL_2(\C)^2$ with $A^2 B^2 = -I$, every trace-zero matrix $C \in \SL_2(\C)$ satisfies $C^2 = -I$ automatically, so $(A,B,C) \in R(\pi_1(N))$ for any such $C$.
\end{proposition}
\begin{proof}
If $\tr(C) = 0$, the characteristic polynomial of $C$ is $t^2 - \tr(C)t + 1 = t^2 + 1$.
By Cayley--Hamilton (Lemma~\ref{lem:traceids}(IV)), $C^2 + I = 0$, i.e.\ $C^2 = -I$.
Substituting into $A^2 B^2 C^2 = I$ gives $A^2 B^2(-I) = I$, i.e.\ $A^2 B^2 = -I$.
Conversely, any $C$ with $\tr(C) = 0$ has characteristic polynomial $t^2 + 1$, so $C^2 = -I$ by Cayley--Hamilton regardless of which trace-zero matrix $C$ is chosen; if $A^2 B^2 = -I$ then $A^2 B^2 C^2 = (-I)(-I) = I$, verifying the surface relation.
\end{proof}

\begin{remark}\label{rem:boundarydim}
Proposition~\ref{prop:boundary} shows that the fiber of $R(\pi_1(N))$ over $\{\tr C = 0\}$ decomposes as a product of two independent pieces:
\begin{itemize}
\item The space of pairs $(A,B) \in \SL_2(\C)^2$ satisfying $A^2 B^2 = -I$ as a matrix equation.
This is a system of $3$ scalar equations on the $6$-dimensional space of pairs, so the expected dimension of the solution set is $3$.
\item The full trace-zero locus $\{\tr C = 0\} \subset \SL_2(\C)$, which is itself $2$-dimensional ($\SL_2(\C)$ is a $3$-dimensional affine variety over $\C$, and $\tr = 0$ is one condition).
\end{itemize}
We emphasize that the condition $A^2 B^2 = -I$ is strictly stronger than its trace-level consequence $\tr(A^2 B^2) = -2$.
Indeed, $\tr(A^2 B^2) = -2$ is equivalent, via Lemma~\ref{lem:A2B2}, to $abp - a^2 - b^2 + 2 = -2$, i.e.\ $abp - a^2 - b^2 + 4 = 0$, which is precisely the condition $c^2 = 0$ from Theorem~\ref{thm:surfacerel}.
Thus the trace-level equation~\eqref{eq:surface} correctly detects the $c = 0$ slice of $S$ via the algebraic condition $c = 0$, but the trace-level data alone is insufficient to determine which pairs $(A,B)$ with $\tr(A^2 B^2) = -2$ actually satisfy the matrix identity $A^2 B^2 = -I$.
A complete treatment of this locus would require a full analysis of the matrix equation $A^2 B^2 = -I$ and the GIT quotient by the $\SL_2(\C)$-action.
We leave this classification as an open problem.
\end{remark}


\end{document}